\documentclass[12pt]{article}
\usepackage{amsfonts}
\usepackage{amsmath,amssymb, amsthm}
\usepackage{euscript}
\usepackage{enumerate}
\usepackage{verbatim}
\usepackage[latin1]{inputenc}

\newtheorem{thm}{Theorem}[section]
\newtheorem{thm*}{Theorem }
\newtheorem{prop}[thm]{Proposition}
\newtheorem{lem}[thm]{Lemma}
\newtheorem{Cor}[thm]{Corollary}

\newtheorem{rem}[thm]{Remark}

\def\bee{\begin{eqnarray}}
\def\bes{\begin{eqnarray*}}
\def\eee{\end{eqnarray}}
\def\ees{\end{eqnarray*}}
\def\a{\alpha}
\def\b{\beta}

\def\s{\sigma}
\def\f{\varphi}
\def\Z{\mathbf Z}
\def\la{\langle}
\def\ra{\rangle}

\def\Proof{{\it Proof.\ }\ }
\def\ctd{\hfill$\Box$}

\title{Irreductible Representations of the Simple Jordan Superalgebra of Grassmann Poisson bracket}

\author{Olmer Folleco Solarte\thanks{Supported by CAPES grant and CNPq grant}\\
\it Universidad del Cauca,\\
\it Calle 5, no.4 - Popay\'an - Colombia\\
{\it E-mail\/}: olmer@unicauca.co\\
{}\\
and\\
{}\\ and Ivan Shestakov\thanks{Supported by FAPESP grant and CNPq grant}\\
 \it Instituto de Matem\'atica e Estat\'istica,\\
\it Universidade de S\~ao Paulo,\\
\it Rua do Mat\~ao, 1010 - CEP 05508-090 - S\~ao Paulo - SP - Brasil\\[1mm]
and\\[1mm]
\it Sobolev Institute of Mathematics,\\
\it Novosibirsk - 630090 - Russia\\
{\it E-mail\/}: shestak@ime.usp.br}
\date{}

\begin{document}

\maketitle

\abstract{We obtain classification of the  irreducible bimodules over the Jordan superalgebra $Kan(n)$, the Kantor double of the Grassmann Poisson superalgebra $G_n$ on $n$ odd generators, for all $n \geq 2$ and an algebraically closed field of characteristic $\ne 2$. This generalizes the corresponding result of C.Mart\'\i nez and E.Zelmanov announced in \cite{MZ2} for $n>4$ and a field of characteristic zero.} 



\section{Introduction}
\hspace{\parindent}
An algebra  $J$ over a field $F$ of characteristic $\neq 2$ is called a {\em Jordan algebra} if it satisfies the identities
\bes
xy &=& yx\\
(x^2 y)x &=& x^2 (yx).
\ees
These algebras were introduced in \cite{JNW} as an algebraic formalism of quantum mechanics. Since then, they have found various applications in mathematics and theoretical physics and now form an intrinsic part of modern algebra. We refer the reader to the books \cite{Jac, Mc, ZSSS} and the survey \cite{KS} for principal results on the structure theory and representations of Jordan algebras.

\smallskip
Jordan superalgebras appeared in 1977--1980 \cite{Kap, Kac}. A {\em Jordan superalgebra} is a $\mathbb{Z}_2$-graded algebra $J = J_{{0}} + J_{{1}}$ satisfiying the graded identities:
\bee   
&xy = (-1)^{|x| |y|} yx,&\nonumber\\
&((xy)z)t + (-1)^{|y||z|+|y||t|+|z||t|}((xt)z)y + (-1)^{|x||y|+|x||z|+|x||t|+|z||t|}((yt)z)x =& \nonumber\\
&=(xy)(zt) + (-1)^{|y||z|}(xz)(yt) + (-1)^{|t|(|y|+|z|)}(xt)(yz),&\label{J}
\eee
where $|x|=i$ if $x \in J_i$. The subspaces $J_{0}$ and $J_{1}$ are referred as the even and the odd parts of $J$, respectively. The even part $J_{ 0}$ is a Jordan subalgebra of $J$, and the odd part $J_{ 1}$ is a Jordan bimodule over $J_{ 0}$. 

\smallskip
In \cite{Kac} (see also \cite{Kan}),  the simple finite dimentional Jordan superalgebras over an algebraically closed field of zero characteristic were classified. The only superalgebras  in this classification whose even part is not semisimple are the  Jordan superalgebras of Grassmann Poisson brackets $Kan(n)$, defined below.

\smallskip
A {\em dot-bracket superalgebra} $A=(A_0 + A_1,\cdot,\{ ,\} )$ is an associative supercommutative superalgebra $(A,\cdot)$ together with a super-skew-symmetric bilinear product $ \{ ,\} $. One can  constract the {\em Kantor superalgebra} $J(A)$ via the {\em Kantor doubling process} as follows \cite{KMc}: Consider the vector space direct  sum  $J=A\oplus \overline{A}$, where $\overline{A}$ is just $A$ labelled, multiplication in $J(A)$ is given by:
\begin{center}
$f\bullet g = f \cdot g$,\\
$f\bullet \overline{g}=\overline{f \cdot g}$,\\
$\overline{f}\bullet g=(-1)^{|g|}\overline{f \cdot g}$,\\
$\overline{f}\bullet \overline{g}=(-1)^{|g|}\{ f ,g \} $,
\end{center}
for $f,g \in A_0 \cup A_1$. Then, we have the $\mathbb{Z}_2$-grading $J(A)=J_0 + J_1$, where $J_0=A_0+\overline{A_1}$ and $J_1=A_1+\overline{A_0}$,  and $J(A)$ is a supercommutative superalgebra  under this grading.
\smallskip

{\bf Theorem (\cite{KMc1})}. {\em If $A = A_0 + A_1$ is a unital dot-bracket superalgebra, then $J(A)$ is a Jordan superalgebra if and only if the following identities hold:
\bee\label{JS}
&\{ f,(g \cdot h)\}=\{ f,g\}\cdot h + (-1)^{|f||g|}g\cdot \{ f,h\} - D(f)\cdot g\cdot h,\label{gleibniz}\\
&\{ f,\{ g,h\}\} - \{ \{ f,g\},h\} - (-1)^{|f||g|}\{ g,\{ f,h\}\} =&\nonumber\\
& D(f)\cdot \{ g,h\} + (-1)^{|g|(|f|+|h|)}D(g)\cdot \{ h,f\} + (-1)^{|h|(|f|+|g|)}D(h)\cdot \{ f,g\}&\label{gjacob}\\
&\{\{ x,x \},x \} = -\{ x,x \}\cdot D(x)&,
\eee
where $D(f)=\{ f,1 \}$, $f,g,h \in A_0 \cup A_1$ and $x \in A_1$. The last identity is needed only in characteristic $3$ case.
}

\smallskip

A dot-bracket  superalgebra $P$ is called a {\em Poisson superalgebra} if it satisfies the identities of the above theorem with $D\equiv 0$. The above construction was first introduced by I.Kantor \cite{Kan} for the Grassmann Poisson superalgebras.

\smallskip

Let $G_n$ be the {\it Grassman superalgebra} generated by $n \geq 2$ odd generators $e_1,e_2,\ldots ,e_n$ over a field $F$, such that $e_ie_j + e_je_i=0$ and $e_i^2=0$. Define on $G_n$ an odd superderivation $\dfrac{\partial}{\partial e_j}$ for $j=1,2,\dots ,n$ by the equalities
\bes
\dfrac{\partial e_i}{\partial e_j} &=& \delta_{ij},\\
\dfrac{\partial(uv)}{\partial e_j}&=&\dfrac{\partial u}{\partial e_j}v + (-1)^{|u|}u\dfrac{\partial v}{\partial e_j},
\ees
and then define the superbracket
\begin{center}
$\{ f,g \}=(-1)^{|f|}\sum\limits_{i=1}^{n}\dfrac{\partial f}{\partial e_i}\dfrac{\partial g}{\partial e_j}$.
\end{center}

One can check that the dot-bracket superalgebra $G_n$ is a Poisson superalgebra, and it was proved in \cite{Kan} that the superalgebra $Kan(n)=J(G_n)$ is a simple Jordan superalgebra for all $n\geq 2$.

\smallskip

A  {\em Jordan (super)bimodule} over  a Jordan superalgebra $J$ is defined in a usual way: a $\mathbb{Z}_2$-graded $J$-bimodule $V=V_0+V_1$  is called a  Jordan $J$-bimodule if the {\em split null extention}  $E(J,V)=J\oplus V$ is a Jordan superalgebra. 
Recall that the multiplication in the split null extention extends the multiplication in $J$ and the action of $J$ on $V$, while the product of two arbitrary elements in $V$ is zero.

The first main problem of the representation theory for any class of algebras is the classification of irreducible bimodules. The description of unital irreducible  finite dimensional Jordan bimodules is  practically finished for  simple finite dimensional Jordan superalgebras over an algebraically closed field of characteristic zero \cite{MS, MSZ, MZ1, MZ2, MZ3, T1, T2, Sht1, Sht2}.   One of  main tools used in these papers was the famous Tits-Kantor-Koecher functor ($TKK$-functor) which associates with a Jordan (super)algebra $J$ a certain Lie (super)algebra $TKK(J)$. Using the known classification of irreducible Lie bimodules over $TKK(J)$, the authors recovered the structure of irreducible bimodules over $J$. Observe that this method may be used only in the characteristic zero case since classification of irreducible Lie supermodules is not known for positive characteristic.  

The classification of irreducible bimodules over the Kantor superalgebra $Kan(n)$ was first obtained in \cite{Sht1} via $TKK$-functor. In \cite{MZ2}, the authors pointed out that the using of the $TKK$-functor in \cite{Sht1} was not quite correct, and the classification obtained there was not complete. They announced a new classification of irreducible bimodules over $Kan(n)$ for all $n>4$ and an algebraically closed field $F$ of characteristic zero.

\smallskip
 In this paper, we classify the irreducible bimodules for the superalgebra $Kan(n)$ over an algebraically closed field $F$ of characteristic $\ne 2$ and $n\geq 2$. 
Our proof is direct and does not use the structure of Lie modules over the Lie superalgebra $L=TKK(Kan(n))$. In order to prove that the constructed bimodule is Jordan, we give a new construction of a Jordan bracket on the tensor product of a  Poisson superalgebra $P$ with a certain generalized derivation and an associative commutative algebra with a derivation.

\section{Some Properties}
\hspace{\parindent}

Recall that the Grassmann algebra $G_n$ has a base formed by 1 and the products $e_{i_1}e_{i_2}\cdots e_{i_n}$, where $1\leq i_1<i_2<\cdots<i_n\leq n$. 

For an ordered subset $I=\{ i_1,i_2,\ldots , i_k \}\subseteq I_n = \{ 1,2,\ldots , n\}$,
we denote 
\bes
e_{I}:=e_{i_1}e_{i_2}\cdots e_{i_k}, 
\ees
so
\bes
\overline{e_I}=\overline{e_{i_1}e_{i_2}\cdots e_{i_k}}, \ \ e_{\phi}=1, \hbox{ and }\overline{e_{\phi}}=\overline{1}.
\ees
Now, as $e_ie_j = -e_je_i$, for $i,j \in I_n$, $i \ne j$, if $\sigma$ is a permutation of the set $I$, we have
\begin{center}
$e_I = sgn(\sigma) e_{\sigma(I)}$,
\end{center}
where $sgn(\sigma)$ is the sign of the permutation $\sigma$.\\

For ordered subsets $I=\{i_1,\ldots i_k \}$ and $J=\{j_1,\ldots j_s \}$, denote by $I\cup J$ the ordered set  
\bes
I\cup J=\{i_1,\ldots i_k,j_1,\ldots j_s \}.
\ees
Then the multiplication in $Kan(n)$ is given by:
\begin{center}
\[
e_{I}\bullet e_{J}=e_{I}e_{J}=
\left\{
\begin{tabular}{ccl}
		$e_{I\cup J}$& if & $I\cap J = \phi, $\\
		$0$& if & $I\cap J \ne \phi, $\\	
\end{tabular}\right.
\]
\\
\[
e_{I}\bullet \overline{e_{J}}=\overline{e_{I}e_{J}}=
\left\{
\begin{tabular}{ccl}
		$\overline{e_{I\cup J}}$& if & $I\cap J = \phi, $\\
		$0$& if & $I\cap J \ne \phi, $\\	
\end{tabular}\right.
\]
\\
\[
\overline{e_{I}}\bullet e_{J}=(-1)^{s}\overline{e_{I}e_{J}}=
\left\{
\begin{tabular}{ccl}
		$(-1)^{s}\overline{e_{I\cup J}}$& if & $I\cap J = \phi, $\\
		$0$& if & $I\cap J \ne \phi, $\\	
\end{tabular}\right.
\]
\\
\[
\overline{e_{I}}\bullet \overline{e_{J}}=(-1)^{s}\{ e_{I},e_{J} \}=
\left\{
\begin{tabular}{ccl}
		$(-1)^{s+k+p+q}e_{I'\cup J'}$& if & $I\cap J = \{ i_p\} =\{ j_q \},$\\
		$0$&  & otherwise,\\	
\end{tabular}\right.
\]
\end{center}
where $I'=\{ i_1,\ldots, i_{p-1}, i_{p+1} ,\ldots, i_k \}$ and $J'=\{ j_1 ,\ldots, j_{q-1}, j_{q+1} ,\ldots, j_s \}$. 

We will use the notation $\bullet$  only in the presence of other multiplications.

\smallskip

Let $V$ be a Jordan bimodule over $Kan(n)$. For $a\in Kan(n)$ we denote by $R_a$ the action of $a$ on $V$:  $R_a(v)=v\cdot a$. 
The Jordan superidentity (\ref{J}) implies the following operator relations:
\bee
&&R_yR_zR_t+(-1)^{|y||z|+|y||t|+|z||t|}R_tR_zR_y+(-1)^{|z||t|}R_{(yt)z}\nonumber\\
&=&R_yR_{zt}+(-1)^{|y||z|}R_zR_{yt}+(-1)^{|t||yz|}R_tR_{yz},\label{RJ1} \\
&&[R_{xy},R_{z}]_s+(-1)^{|y||z|}[R_{xz},R_{y}]_s+(-1)^{|x||yz|}[R_{yz},R_{x}]_s=0, \label{RJ2}
\eee
where $[R_{x},R_{y}]_s = R_{x}R_{y}-(-1)^{|x||y|}R_{y}R_{x}$ denotes the supercomutator of the operators $R_x,R_y$.

\smallskip

It is well known (see, for instance, \cite{Jac, MZ1}) that every Jordan bimodule $V$ over a unital Jordan (super)algebra $J$ is decomposed into a direct sum of three subbimodules
\[
V=V(0)\oplus V(1)\oplus V(1/2),
\]
where $V(0)$ is a trivial bimodule, $V(1)$ is a unital bimodule, and $V(1/2)$ is a {\em semi-unital} or {\em one-sided} bimodule, that is, where the unit 1 of $J$ acts as $\tfrac12$. Moreover, for a semi-unital bimodule $V$, the mapping $a\mapsto 2R_a$ is a homomorphism of a Jordan (super)algebra $J$ into the special Jordan (super)algebra $(End\,V)^{+}$. Therefore, a simple exceptional  unital Jordan (super)algebra admits only unital irreducible bimodules.

It was shown in \cite{Sh} that the Kantor double $J(P)$ for a Poisson superalgebra $P$ is special if and only if it satisfies the identity $\{\{P,P\},P\}=0$. Evidently, the superalgebra $G_n$ does not satisfy this condition, hence the superalgebra $Kan(n)=J(G)$ is exceptional. In particular, every irreducible Jordan bimodule over $Kan(n)$ is unital.

Below $V$ denotes a unital Jordan bimodule over the superalgebra $Kan(n)$. 

\smallskip

The next Lemma gives the first properies of the right operators over $V$.

\begin{lem}\label{lem1}
Given  index sets $I=\{ i_1,\ldots ,i_k \}$ and $J=\{ j_1,\ldots ,j_s \}$  contained in $I_n = \{ 1,\ldots , n \}$, we have
\begin{enumerate}
\item $[R_{e_I},R_{e_J}]_s =0$, for all $I$ and $J$.
\item $[R_{e_I},R_{\overline{e_J}}]_s = 0$, if $|J\cap I|\geq 2$.
\item $[R_{e_I},R_{\overline{1}}]_s = 0$, for all $I\ne \{ 1,2,\ldots , n \}$.
\item $[R_{\overline{e_I}},R_{\overline{e_J}}]_s =0$, if $I\cap J \ne \phi$.
\end{enumerate}
\end{lem}

\Proof
We will leave the proof of item $1$ to the end.

For item $2$, the set $I\cap J$ has at least two elements, and we can assume, without loss of generality, that these elements are $j_1$ and $j_s$. Let
\bes
x = e_{j_1}, \ y = \overline{e_{j_2}\cdots e_{j_s}},  \hbox{ and } z = e_I,
\ees
then
\begin{center}
$xy = e_J$ and $xz = yz = 0$,
\end{center}
and relation (\ref{RJ2}) finishes the proof.

For item $3$, as $I\ne \{1,2,\ldots , n \}$, there is $p\notin I$, so if
\begin{center}
$x = \overline{e_{p}}$,\ $\ y = \overline{e_{p}e_{I}}$ and $z = \overline{1}$,
\end{center}
we have
\begin{center}
$xy = (-1)^{k}e_I$ and $xz = yz = 0$,
\end{center}
and again by (\ref{RJ2}) the item is shown.

For item $4$, we take
\begin{center}
$x = e_{I}$, $\ y= \overline{1}$ and $z = \overline{e_{J}}$,
\end{center}
then
\begin{center}
$xy = \overline{e_{I}}$ and $xz = yz = 0$,
\end{center}
and one more time using (\ref{RJ2}) we obtain the result.

For item $1$, first we suppose that $e_{I}\ne e_{J}$, so there exists $e_{j_1}$ such that $e_{j_1}\notin I$. Therefore,  taking
\begin{center}
$x=\overline{e_{j_1}}$, $\ y=\overline{e_{j_1}e_I}$ and $z=e_{J}$,
\end{center}
we have
\begin{center}
$xy=(-1)^{k}e_I$ and $xz = yz = 0$, 
\end{center}
and (\ref{RJ2}) proves the item.

Now, if $I=J=\{ e_i \}$, we take $e_j \ne e_i$, and
\begin{center}
$x = \overline{e_j}$, $\  y= \overline{e_j e_i}\ $ and $z = e_{i}$,
\end{center}
then
\begin{center}
$xy = e_i$, $\ xz = -\overline{e_j e_i}$ and $yz = 0$,
\end{center}
hence by (\ref{RJ2}) and  item $4$, 
\bes
[R_{e_i},R_{e_i}]_s= [R_{xy},R_z]_s=[R_y,R_{xz}]_s=[R_{\overline{e_j e_i}}, R_{\overline{e_j e_i}}]_s = 0.
\ees
Finally, if $I=J$ and $|I|\geq 2$, we take
\begin{center}
$x=e_{i_1}$, $y=e_{i_2}\cdots e_{i_k}$, and $z=e_{I}$,
\end{center}
then
\begin{center}
$xy = e_{I},\ xz = yz = 0$,
\end{center}
and  (\ref{RJ2}) finishes the proof.
\ctd

As a corollary, we obtain the next lemma.

\begin{lem}\label{lem2}
The following statements hold:
\begin{enumerate}
\item If $a \in Kan(n)_1$, $a = e_I$ or $\overline{e_I}$ and $a \ne \overline{1}$, then
$R^2_a = 0$.
\item If $a \in Kan(n)_0$, $a = e_I$ or $\overline{e_I}$ and $a \ne 1,\ \overline{e_i}$, then $R^3_a = 0$.
\item $R^3_{\overline{e_i}} = R_{\overline{e_i}}$, for all $i \in \{ 1,\ldots , n \}$.
\item If $V$ is irreducible and $F$ is algebrically closed then
$R^2_{\overline{1}} = \alpha$ for some $\alpha \in F$.
\end{enumerate}
\end{lem}
\Proof
By items $1$ and $4$ of the previous lemma, for $a \in Kan(n)_1,\ a = e_I$ or $a=\overline{e_I}$,  and $a \ne \overline{1}$ we have
\begin{center}
$[R_a,R_a]_s =2R_a^2= 0$,
\end{center}
which proves  item $1$.

Now, if $a \in Kan(n)_0$,  by superidentity (\ref{RJ1}) we have
\bes
2R_a^3+R_{a^3}-3R_aR_{a^2}=0.
\ees
If $a = e_I$ or $\overline{e_I}$ and $\ a \ne 1,\ \overline{e_{i}}$, then $a^2=0$ and  $R^3_a = 0$, proving  item $2$. 

On the other hand,  if $a = \overline{e_i}$ then $a^2=1$, and since $V$ is unital, the same identity implies 
$2R^3_a=2R_a$, proving item 3.

For item $4$, we recall  the following identity which holds in Jordan algebras \cite{Jac}:
\bes
(a,d,b)c-(a,dc,b)+d(a,b,c) = 0,
\ees
where $(a,b,c)=(ab)c-a(bc)$ is the associator of the elements $a,b,c$.
For Jordan superalgebras, the super-version of this identity holds:
\bes
(a,d,b)c-(-1)^{|b||c|}(a,dc,b)+(-1)^{|a||d|}d(a,b,c) = 0.
\ees
Now, if we take $c \in Kan(n),\ a = b = \overline{1}$, and $d=v\in V$, is easy to see that $(a,c,b)=0$, hence
\begin{center}
$(\overline{1},v,\overline{1})c = (-1)^{|c|}(\overline{1},vc,\overline{1})$.
\end{center}
Therefore, $U = (\overline{1},V,\overline{1})$ is a subbimodule of $V$, and as $V$ is irreducible, we have $U = 0$ or $U = V$.

If $U = 0$, it is clear that $R^2_{\overline{1}} = 0 \in F$. Otherwise $R^2_{\overline{1}}$ is an authomorphim of $V$, and by the Schur lemma, $R^2_{\overline{1}} = \alpha \in F$.
 \ctd

\section{Special Element in $V$}

\begin{lem}\label{lem3}
If $V$ is an unital Jordan bimodule over $Kan(n)$, then there exists $0 \ne v \in V_0\cup V_1$ such that
\begin{center}
$ve_{I} = v\overline{e_I} = 0$,
\end{center}
for all $\phi \ne I \subseteq I_n = \{ 1,\ldots , n  \}$.
\end{lem}
\Proof
For $w \in V$, denote $N_w = \{ a \in Kan(n) \,|\, wa=0 \}$.
We want to find $0\ne v\in V$ such that $e_I,\overline{e_I}\in N_v$ for all $\phi \ne I \subseteq I_n $.

\smallskip

As $[R_{e_I},R_{e_J}]_s = 0$ for all $I,J\subseteq I_n$, and $R^3_{e_I}=0$ for all $I\neq\phi$, the subsuperalgebra of $End\, V$ generated by the set $\{R_{e_I}\,|\,\phi \ne I\subseteq I_n\}$ is nilpotent. Therefore,  there exists $0 \ne u \in V_0\cup V_1$ such that
\begin{center}
$e_I \in N_u$, for all $\phi \ne I\subseteq I_n$.
\end{center}

If $\overline{e_{I_n}}\notin N_u$, consider $u_1 = u\overline{e_{I_n}}$. Since $[R_{e_I},R_{\overline{e_{I_n}}}]_s = 0$ for $|I|\geq 2$,  for these $I$'s we have  $e_I \in N_{u_1}$. In order to show that $e_i \in N_{u_1}$ for all $i \in I_n$, we first substitute in the main Jordan superidentity (\ref{J}) 
$x=\overline{e_i}$, $y = e_{I'_n}$, $z = u$, and \ $t = e_i$, where $I'_n = \{ 1,\ldots , i-1,i+1,\ldots , n \}$.
Then we obtain
\begin{center}
$(u\overline{e_{I_n}})e_i = (u\overline{e_i})e_{I_n}$.
\end{center}
Substituting now again in (\ref{J})
$x = u$, $y = e_{I'_n}$, $z = \overline{e_i}$, and \ $t = e_i$,
 we get
\begin{center}
$(u\overline{e_i})e_{I_n} = 0$,
\end{center}
hence $u_1 e_i = 0$, for all $i \in I_n$.
Therefore,
\begin{center}
$e_I \in N_{u_1}$,  \ for all $\phi \ne I \subseteq I_n$.
\end{center}
If $\overline{e_{I_n}}\notin N_{u_1}$,  we consider the element $u_2=u_1\overline{e_{I_n}}$ and again get  
\begin{center}
$e_I \in N_{u_2}$, for all $\phi \ne I \subseteq I_n$.
\end{center}
Since $R^{3}_{\overline{e_{I_n}}} = 0$, we conclude that there exists $0\neq w\in \{u,u_1,u_2\}$ such that 
\bes
e_J,\ \overline{e_{I_n}} \in N_w \hbox{ for all } \phi \ne I \subseteq I_n.
\ees
For elements $\overline{e_I}$ with $2\leq |I| < n$, substitute in (\ref{J})
$x = e_{i_1}$, $y = e_{I'}$, $z = \overline{1}$,  \ and $t = w$,\ 
 where $I = \{ i_1,\ldots , i_k\}$ and $I'= \{ i_2,\ldots , i_k\}$; then we get
\begin{center}
$w \overline{e_I} = (-1)^{k+1}(w\overline{1})e_I$.
\end{center}
Since 
$[R_{\overline{1}},R_{e_I}]_s = 0$ and $w e_I =0$,
this implies
\begin{center}
$\overline{e_I} \in N_{w}$ for all $I$ with $|I|\geq 2$.
\end{center}

At this point, we need to incorporate the elements  $\overline{e_i}$ \ for $i \in I_n$. First we show that
\begin{center}
$(w\overline{e_i})\overline{e_i} = (w \overline{e_j})\overline{e_j}$, \ for all $i \ne j$.
\end{center}
Substituting in (\ref{J}) 
$x = w$, $y = \overline{e_i}$, $z = \overline{1}$, and $t = \overline{e_j}$,\ 
we obtain
\begin{center}
$((w \overline{e_i})\overline{1})\overline{e_j} = - ((w \overline{e_j})\overline{1})\overline{e_i}$,
\end{center}
and continuing with
$x = w \overline{e_i}$, $y = \overline{1}$, $z = \overline{e_j}$, and $t = \overline{e_i e _j}$, \ we get
\begin{center}
$(((w \overline{e_i})\overline{1})\overline{e_j})\overline{e_i e_j}  - (((w \overline{e_i})\overline{e_i e_j})\overline{e_j})\overline{1} = ((w \overline{e_i})\overline{1})e_i$.
\end{center}
Since $[R_{\overline{e_i}}, R_{\overline{e_ie_j}}]_s = 0$ we have
\begin{center}
$(((w \overline{e_i})\overline{e_i e_j})\overline{e_j})\overline{1} = (((w \overline{e_i e_j})\overline{e_i})\overline{e_j})\overline{1} = 0$,
\end{center}
so
\begin{center}
$(((w \overline{e_i})\overline{1})\overline{e_j})\overline{e_i e_j} = ((w \overline{e_i})\overline{1})e_i$.
\end{center}
Substituting again in (\ref{J})
$x = w$, $y = \overline{e_i}$, $z = \overline{1}$, and $t = e_i$,
we obtain
\bee\label{wee}
((w \overline{e_i})\overline{1})e_i = -(w \overline{e_i})\overline{e_i},
\eee
hence
\begin{center}
$(((w \overline{e_i})\overline{1})\overline{e_j})\overline{e_i e_j} = -(w \overline{e_i})\overline{e_i}$.
\end{center}
Similarly,
\begin{center}
$(((w \overline{e_j})\overline{1})\overline{e_i})\overline{e_i e_j} = -((w \overline{e_j})\overline{1})e_j = (w \overline{e_j})\overline{e_j}$,
\end{center}
and finally
\bes
-(w \overline{e_i})\overline{e_i} = (w \overline{e_i})\overline{1})\overline{e_j})\overline{e_i e_j}= -(((w \overline{e_j})\overline{1})\overline{e_i})\overline{e_i e_j}= -(w \overline{e_j})\overline{e_j},
\ees
what we wanted to prove.

Now, since $R^3_{\overline{e_i}} = R_{\overline{e_i}}$,  we have
\begin{center}
$((w \overline{e_i})\overline{e_i})\overline{e_j} = ((w\overline{e_j})\overline{e_j})\overline{e_j} = w \overline{e_j}$.
\end{center}
Therefore, if there exists $k \in I_n$ such that $w \overline{e_k} = 0$, then $w \overline{e_i} = 0$ for all $i \in I$, and we finish the proof taking $v = w$.

\smallskip
Suppose then that $w \overline{e_i} \ne 0$ for all $i \in I_n$, then we show that
\begin{center}
$(w \overline{e_i})e_i \ne 0$,\ for all $i \in I_n$,
\end{center}
In fact, by item 3 of Lemma \ref{lem1}, $R_{\bar 1}R_{e_i}=-R_{e_i}R_{\bar 1}$, hence by (\ref{wee})
\[
 (w \overline{e_i})\overline{e_i}=((w \overline{e_i})e_i)\overline{1},
\]
and if $(w \overline{e_i})e_i = 0$ for some $i \in I_n$, then
\begin{center}
$0 \ne w \overline{e_i} = ((w \overline{e_i})\overline{e_i})\overline{e_i} = (((w \overline{e_i})e_i)\overline{1})\overline{e_i} = 0$,
\end{center}
a contradiction. Therefore,
\begin{center}
$(w \overline{e_i})e_i \ne 0$ \ for all $i \in I_n$,
\end{center}
Furthermore, for $i,j \in I_n$ with $i \ne j$ we substitute in (\ref{J})
$x = w$, $y = \overline{e_i}$, $z = \overline{e_i e_j}$, and $t = \overline{e_j}$,
then in view of the relations $[R_{e_i},R_{e_ie_j}]_s = [R_{e_j},R_{e_ie_j}]_s = 0$, we obtain
\begin{center}
$(w \overline{e_i})e_i = (w \overline{e_j})e_j$.
\end{center}

We now show that the element   $v=(w \overline{e_1})e_1$ satisfies the statement of lemma. 
Let $i \in I_n$, then
\begin{center}
$v e_i = ((w \overline{e_1})e_1)e_i = ((w \overline{e_i})e_i)e_i = 0$, since $R^2_{e_i} = 0$.
\end{center}
Now, let $I = \{ i_1,\ldots , i_k \}\subseteq I_n$ be, with $k \geq 2$, then in view of item 1 of Lemma \ref{lem1}
\begin{center}
$v e_I = ((w \overline{e_{i_1}})e_{i_1})e_I = \pm ((w \overline{e_{i_1}})e_I)e_{i_1}$.
\end{center}
Substituting in (\ref{J})
$x = w$, $y = e_{I'}$, $z = \overline{e_{i_1}}$, and $t = e_{i_1}$, where $I' = \{ i_2,\ldots , i_k  \}$,
we obtain
\begin{center}
$(w \overline{e_{i_1}})e_I = 0$,
\end{center}
hence
\begin{center}
$v e_I = 0$, \ for all $\phi \ne I \subseteq I_n$.
\end{center}

Analogously as we showed that
$w \overline{e_I} = 0$, for $I \subseteq I_n$ with $2\leq |I| = k <n$,
we can show that
$v \overline{e_I} = 0$   for these $I$'s.  Furthermore, substituting in (\ref{J})
$x = w$, $y = \overline{e_1}$, $z = e_1$, and \ $t = \overline{e_{I_n}}$,
 we have
\begin{center}
$v \overline{e_{I_n}} = ((w\overline{e_1})e_1)\overline{e_{I_n}} = 0$.
\end{center}
Finally, substituting in (\ref{J})
$x = w$, $y = \overline{e_i}$, $z = e_i$, and $t = \overline{e_i}$, 
we obtain
\begin{center}
$v\overline{e_i} = ((w \overline{e_i})e_i)\overline{e_i} = 0$,
\end{center}
ending the proof.
\ctd

\section{Action of $Kan(n)$ on $V$.}

\hspace{\parindent}
 In this section, we will assune that the bimodule $V$ is irreducible. We will find a finite set that generates $V$ as a vector space and will determine the action of the superalgebra $Kan(n)$ on this set.

Let us begin with notation. If $I = \{ i_1, \ldots , i_k \} \subseteq I_n = \{ 1, \ldots , n \}$ and $w\in V$, we denote
\begin{center}
$w(I) := w\overline{1}\overline{e_{i_1}}\overline{1}\cdots \overline{1}\overline{e_{i_k}} := (\cdots(((w\overline{1})\overline{e_{i_1}})\overline{1})\cdots \overline{1})\overline{e_{i_k}}$,
\end{center}
and
\begin{center}
$\overline{w(I)} := w\overline{1}\overline{e_{i_1}}\overline{1}\cdots \overline{1}\overline{e_{i_k}}\overline{1} := ((\cdots(((w\overline{1})\overline{e_{i_1}})\overline{1})\cdots \overline{1})\overline{e_{i_k}})\overline{1}$.
\end{center}
In particular, $w(\phi) = w$ and $\overline{w(\phi)} = w\overline{1}$.

It follows from (\ref{RJ1}) that
\bee\label{RRR}
R_{\bar e_i}R_{\bar{1}}R_{\bar {e}_i} = 0, 
\eee
and
\bee\label{RRR1}
R_{\bar{e}_i}R_{\bar{1}}R_{\bar{e}_j} = -R_{\bar{e}_j}R_{\bar{1}}R_{\bar{e}_i}, \ \hbox{ for } i \ne j,
\eee
so if $\sigma$ is a permutation of $I$,
\begin{center}
$w(I) = sgn(\sigma) w(\sigma(I))$,
\end{center}
where $sgn(\sigma)$ is the sign of $\sigma$.

\smallskip

We want to show that the subspace of $V$ generated by the elements 
$v(I)$ and $\overline{v(I)}$, where $I$ runs all the subsets of $I_n$ and $v$ is the element from the previous section, is closed under the action of $Kan(n)$ and hence coincides with $V$.

\begin{lem}\label{lem4}
If $I, J \subseteq I_n$, with  $J\not\subseteq I$, then
\begin{center}
$v(I)e_J = v(I)\overline{e_J} = \overline{v(I)}e_J = 0$.
\end{center}
Moreover, if $|J\setminus I|\geq 2$, then
\begin{center}
$\overline{v(I)}\overline{e_J} = 0$.
\end{center}
\end{lem}
\Proof
Let $I = \{ i_1,\ldots , i_k \}$ and $J = \{ j_1,\ldots , j_s \}$.
We use induction on $|I| = k$. If $k = 0$, then by the properties of the element $v$ we have
\begin{center}
$ve_J = v\overline{e_J} =  0$, \ if $|J|\geq 1$.
\end{center}
Moreover, by item 3 of Lemma \ref{lem1},  $[R_{\overline{1}},R_{e_J}]_s=0$ for $J\neq I_n$, hence
\begin{center}
$(v\overline{1})e_J =\pm (ve_J)\overline{1} = 0$ \ for $J\neq I_n$.
\end{center}
For $J = I_n$, we substitute in (\ref{J}) 
$x = e_1$, $y = e_{I'}$, $z = \overline{1}$,\  and $t = v$, where $I' = \{ 2,\ldots,n \}$,  then
$(v\overline{1})e_{I_n} = \pm v\overline{e_{I_n}} = 0.$

Finally, substituting in (\ref{J})
$x = \overline{1}$, $y = e_{j_1}$, $z = v$,\ and $t=\overline{e_{J'}}$, where $J' = \{ j_2,\ldots , j_s \}$,
we obtain
\begin{center}
$(v\overline{1})\overline{e_J} = 0$, if $|J|\geq 2$.
\end{center}
\smallskip

Now, suppose that the lemma is true for $|I|=m < k$. Let $I' = I\setminus\{ i_{k} \}$, then we have by (\ref{J}) and induction on $|I|$
\bes
v(I){e_{J}}&=&((v(I')\bar 1)\bar{e}_{i_k})e_{J}\\
&=&\pm ((v(I')e_{J})\bar{e}_{i_k})\bar 1\pm  v(I')(\bar{e}_J\,\bar{e}_{i_k})+(v(I')\bar 1)(\bar{e}_{i_k}e_J)\pm (v(I')\bar{e}_{i_k})\bar{e}_J\\
&=& \pm v(I') (\overline{e_J}\,\overline{e_{i_k}})\pm \overline{v(I')}(\overline{e_{i_k}}e_J). 
\ees
Consider the two cases. If $i_k\in J$, then $\overline{e_J}\,\overline{e_{i_k}}=\pm e_{J\setminus\{i_k\}},\ \overline{e_{i_k}}e_J=0$.
Since $J\not\subseteq I$, we have $J\setminus\{i_k\}\not\subseteq I'$. Therefore, by induction on $|I|$, ${v(I')} e_{J\setminus\{i_k\}}=0$. Finally, if $i_k\not\in J$ then 
$\overline{e_J}\,\overline{e_{i_k}}=0,\ \overline{e_{i_k}}e_J=\pm \overline{e_{J\cup \{i_k\}}}$. Clearly, 
$|(J\cup\{i_k\})\setminus I'|\geq 2$ hence by induction on $|I|$ we have $\overline{v(I')}\overline{e_{J\cup \{i_k\}}}=0$.   Therefore, in both cases $v(I){e_{J}}=0$, proving  first equality of the lemma.

\smallskip
Similarly,
\bes
v(I)\overline{e_{J}}&=&((v(I')\bar 1)\overline{e_{i_k}})\overline{e_{J}}=\pm ((v(I')\overline{e_{J}})\overline{e_{i_k}})\bar 1\pm (v(I')\bar 1)(\overline{e_{i_k}}\,\overline{e_{J}})\\
&=& \hbox{(by induction on $|I|$)} = \pm\overline {v(I')} (\overline{e_{i_k}}\,\overline{e_{J}}). 
\ees
As in the previous case, we have $\overline {v(I')} (\overline{e_{i_k}}\,\overline{e_{J}})=0$,  proving  second equality of the lemma.

\smallskip

Since $[R_{\overline{1}},R_{e_{J}}]_s = 0$ for $J\ne I_n$, we have $\overline{v(I)}e_{J} = (v(I)\bar 1)e_{J}=\pm (v(I)e_{J})\bar 1=0$, with $J \ne I_n, \ J\not\subseteq I$.
For $J = I_n$ we have by (\ref{J}) 
\bes
\overline{v(I)}e_{I_n}&=& ((\overline{v(I')}\overline{e_{i_k}})\bar 1)e_{I_n}=\pm((\overline{v(I')}e_{I_n})\bar 1)\overline{e_{i_k}}+(-1)^n(\overline{v(I')}\overline{e_{i_k}})\overline{e_{I_n}}\\
&=&\hbox{ (by induction on $|I|$) } =(-1)^nv(I) \overline{e_{I_n}},
\ees
where $I'=I\setminus\{i_k\}$.
Therefore, for any $I\subsetneq I_n$ we have $\overline{v(I)}e_{I_n}=\pm v(I) \overline{e_{I_n}}=0$. We will also need later the following equality for $I=I_n$: 
\bee\label{In}
\overline{v(I_n)}\,e_{I_n}=(-1)^{n} v(I_n) \overline{e_{I_n}}.
\eee
Finally, we have by (\ref{J}) and Lemma \ref{lem4} for $J' = J\setminus\{ j_1 \}$, with $e_{j_1}\not\in I$:
\bes
\overline{v(I)}\,\overline{e_{J}}=(v(I)\bar 1)(e_{j_1}\overline{e_{J'}})=\pm \overline{e_{j_1}}(v(I)\overline{e_{J'}})\pm((\bar 1e_{j_1})v(I))\overline{e_{J'}}\pm \bar 1(v(I)\overline{e_{J}}).
\ees
By the previous cases, since $j_1\not\in I$ and $J'\not\subseteq I$, we have
\bes
&((\bar 1e_{j_1})v(I))\overline{e_{J'}}=\pm (v(I)\overline{e_{j_1}})\overline{e_{J'}}=0,&\\
&v(I)\overline{e_{J}}=v(I)\overline{e_{J'}}=0,&
\ees
hence 
$$
\overline{v(I)}\overline{e_{J}}=0,
$$
proving the lemma.
\ctd

\smallskip
\begin{lem}\label{lem5}
Let $J=\{j_1,\ldots,j_s\} \subseteq I=\{i_1,\ldots,i_{k-s},j_s,j_{s-1},\ldots,j_1\}$. Then 
\begin{itemize}
\item $v(I)e_J = v(I \setminus J)$, 
\item $v(I)\overline{e_J} = \overline{v(I\setminus J)}$,
\item  $\overline{v(I)}e_J = (-1)^{|J|}\overline{v(I \setminus J)}$,
\item $\overline{v(I)}\,\overline{e_J}=(-1)^{|J|-1}\alpha (|J|-1)v(I\setminus J)$,
\end{itemize}
where $\alpha=R_{\bar 1}^2$.
Furthermore, if $|J\setminus I|=1,\ I=\{i_1,\ldots,i_{k-s+1},j_{s-1},\ldots,j_1\},\ J= \{j_1,\ldots,j_s\}$, then 
$$ 
\overline{v(I)}\,\overline{e_J}=(-1)^{s-1}\overline{v(I\setminus J)}\,\overline{e_{j_s}}=(-1)^{s-1}v((I\setminus J)\cup\{j_s\})=(-1)^{s-1}v(\{i_1,\ldots,i_{k-s+1},j_s\}).
$$
\end{lem}
\Proof
We will use induction on $|J| = s$. If $s = 0$, we have $e_J = 1$ and all the claims are clear. For $s = 1$, consider first $v(I)e_{j_1}$. Let $I' =I\setminus\{j_1\}$, then by (\ref{J})  and Lemma \ref{lem4} we have
\bes
v(I)e_{j_1}&=&((v(I')\bar 1)\bar e_{j_1})e_{j_1}\\&=&((v(I')e_{j_1})\bar e_{j_1})\bar 1+v(I')(\bar e_{j_1}\,\bar e_{j_1})-(v(I')\bar e_{j_1})\bar e_{j_1}=v(I'),
\ees
which proves  first equality for $s=1$. Similarly,
\bes
v(I)\overline{e_{j_1}}&=&((v(I')\bar 1)\bar e_{j_1})\bar{e}_{j_1}\\
&=&-((v(I')\bar{e}_{j_1})\bar e_{j_1})\bar 1+\overline{v(I')}(\bar e_{j_1}\,\bar e_{j_1})=\overline{v(I')}.
\ees
Third equality is true for $s=1$ since 
$$
\overline{v(I)}e_{j_1}=(v(I)\bar 1)e_{j_1}=
-(v(I)e_{j_1})\bar 1=-v(I')\bar 1=-\overline{v(I')}.
$$  
Furthermore, it follows from (\ref{RRR}) and (\ref{RRR1}) that $\overline{v(I)}\,\overline{e_{j_1}}=0$, which proves  fourth equality for $s=1$. Finally, if $j\not\in I$ then by definition 
$\overline{v(I)}\,\overline{e_{j}}=v(I\cup\{j\})$, proving the last equality for $s=1$.

\smallskip

Assume now that the lemma is true if $|J|<s$. Let $I' = I\setminus\{j_1 \}$ and $J' = J\setminus\{ j_1 \}$; then by (\ref{J}), Lemma \ref{lem4}, and induction on $|J|$, we have
\bes
v(I)e_J&=&((v(I')\bar 1)\overline{e_{j_1}})e_{J}\\
&=&\pm((v(I')e_J)\overline{e_{j_1}})\bar 1-(-1)^sv(I')(\overline{e_J}\,\overline{e_{j_1}})\pm ( v(I')\overline{e_{j_1}})\overline{e_J}\\
&=&v(I')e_{J'}=v(I'\setminus J')=v(I\setminus J).
\ees
Similarly,
\bes
v(I)\overline{e_J}&=&((v(I')\bar 1)\overline{e_{j_1}})\overline{e_{J}}\\
&=&\pm((v(I')\overline{e_J})\overline{e_{j_1}})\bar 1+(v(I')\bar 1)(\overline{e_J}\,\overline{e_{j_1}})\\
&=&(-1)^{s-1}\overline{v(I')}e_{J'}=\overline{v(I'\setminus J')}=\overline{v(I\setminus J)}.
\ees
Furthermore, if $J\neq I_n$ then
\bes
&\overline{v(I)}e_J=(v(I)\bar 1)e_J=(-1)^{|J|}(v(I)e_J)\bar 1=(-1)^{|J|}v(I\setminus J)\bar 1=(-1)^{|J|}\overline{v(I\setminus J)}, &
\ees
 and for $J=I_n$ we have by (\ref{In}) $\overline{v(I_n)}e_{I_n}=(-1)^nv(I_n)\overline{e_{I_n}}=(-1)^n(v\bar 1)$, which proves  third equality.

To prove fourth equality, observe first that
\bes
v(I)\overline{e_{J'}}&=&v(\{i_1,\ldots,i_{k-s},j_s,\ldots,j_2,j_1\})\overline{e_{\{j_2,\ldots,j_s\}}}\\
&=&(-1)^{s-1}v(\{i_1,\ldots,i_{k-s},j_1,j_s,\ldots,j_2\})\overline{e_{\{j_2,\ldots,j_s\}}}\\
&=&(-1)^{s-1}\overline{v(\{i_1,\ldots,i_{k-s},j_1\})}=(-1)^{s-1}\overline{v(I\setminus J')}.
\ees
Now, applying again (\ref{J}), Lemma \ref{lem4}, and induction on $|J|$, we get
\bes
\overline{v(I)}\overline{e_J}&=& (v(I)\bar 1)({e_{j_1}}\overline{e_{J'}})\\
&=&-(v(I)\overline{e_{J'}})(\bar 1{e_{j_1}})+((v(I)\bar 1)e_{j_1})\overline{e_{J'}}
-((v(I)\overline{e_{J'}})e_{j_1})\bar 1\\
&=&\overline{v(I\setminus J')}\overline{e_{j_1}}+(\overline{v(I)}e_{j_1})\overline{e_{J'}}-(-1)^{s-1}(\overline{v(I\setminus J')}e_{j_1})\bar 1\\
&=&-\overline{v(I')}\overline{e_{J'}}+(-1)^{s-1}\overline{v(I\setminus J)}\bar1\\
&=&-(-1)^{s-2}\alpha (s-2)v(I'\setminus J')+(-1)^{s-1}\alpha v(I\setminus J)\\
&=&(-1)^{s-1}\alpha(s-1) v(I\setminus J).
\ees

Finally, let $I=\{i_1,\ldots,i_{k-s+1},j_{s-1},\ldots,j_1\},\ J= \{j_1,\ldots,j_s\},\ J'=J\setminus\{j_s\}$. Arguing as above, we get by Lemma \ref{lem4}
\bes
\overline{v(I)}\overline{e_J}&=& (v(I)\bar 1)({e_{J'}}\overline{e_{j_s}})\\
&=&-(v(I)\overline{e_{j_s}})(\bar 1{e_{J'}})+((v(I)\bar 1)e_{J'})\overline{e_{j_s}}
-((v(I)\overline{e_{j_s}})e_{J'})\bar 1\\
&=&(\overline{v(I)}e_{J'})\overline{e_{j_s}}=(-1)^{s-1}\overline{v(I\setminus J')}\overline{e_{j_s}}=(-1)^{s-1}v((I\setminus J')\cup\{j_s\}).
\ees
This finishes the proof of the lemma.
\ctd

\smallskip
Lemmas 3 -- 5 imply the next theorem:
\begin{thm}\label{thm1}
Let $V$ be a unital irreducible bimodule over the superalgebra $Kan(n)$ and $v\in V$ be a special element from Lemma \ref{lem3}, then $V$ is generated as a vector space by elements of the type
\bee\label{base}
v(I), \ \overline{v(I)}, \hbox{\ where } I \subseteq I_n = \{ 1,\ldots , n \}.
\eee
 Furthermore, let $I,J\subseteq I_n$,\ $J=\{j_1,\ldots,j_{s_1},j_{s_1+1},\ldots,j_{s_1+s_2}\}$,\  $I=\{i_1,\ldots, i_{k-s_1},j_{s_1},\cdots,j_{1}\}$. Then the action of $Kan(n)$ on $V$ is defined, up to permutations of the index sets $I$ and $J$, as follows:
\bes
v(I) e_{J} &=&
\begin{cases}
		v(I\setminus J) \hbox{\quad if }& s_2=0 \ (\hbox{or, equivalently, } J\subseteq I),\\
		0  & \hbox{ otherwise; }	
\end{cases}\\[1mm]
v(I) \overline{e_{J}} &=&
\begin{cases}
		\overline{v(I\setminus J)} \hbox{\quad if }& s_2=0,\\
		0  & \hbox{ otherwise; }	
	\end{cases}\\[1mm]
\overline{v(I)} e_{J} &=&
\begin{cases}
		(-1)^{s}\overline{v(I\setminus J)}\hbox{ \quad if }& s_2=0,\\
			0  & \hbox{ otherwise; }		
\end{cases}\\[1mm]\overline{v(I)} \overline{e_{J}} &=&
\begin{cases}
		(-1)^{s_1}\overline{v(I\setminus J_1)}\,\overline{e_{J_2}} \hbox{\quad\quad\quad if }& s_2 = 1,\\
		(-1)^{s-1}\alpha (s-1)\,v(I\setminus J)\ \hbox{ if } & s_2 = 0,\\
		0  & \hbox{ otherwise, } \\
\end{cases}
\ees
where $\alpha = R^2_{\overline{1}}$,   and $s = s_1 + s_2 = |J|$.
\end{thm}
\Proof Let $W$ be the vector subspace of $V$ spanned by the set (\ref{base}). It follows from Lemmas \ref{lem4} and \ref{lem5} that $W Kan(n)\subseteq W$, that is, $W$ is a subbimodule of $V$. Clearly, $W\neq 0$, hence $W=V$. The rest of the theorem follows directly from Lemmas \ref{lem4} and \ref{lem5}.

\ctd

Since the action of $Kan(n)$ on $V$ depends on the choise of a special element $v\in V$ and a parameter $\alpha = R_{1}^2\in F$,  we will denote the bimodule $V$ as $V(v,\alpha)$.

\section{Linear independence, Irreducibility, and Isomorphism problem}

\hspace{\parindent}
In this section, we will prove that the set (\ref{base}), when $I$ runs all different (non-ordered) subsets of $I_n$, is in fact linearly independent and hence forms a base of the bimodule $V$. Furthermore, we will prove that $V(v,\alpha)$ is irreducible for any $\a\in F$ and that the bimodules $V(v,\a)$ and $V(v',\a')$ are isomorphic if and only if $|v|=|v'|,\ \a=\a'$.

\begin{lem}\label{lem7}
Given $I \subseteq I_n $, there exists an element $W=W(I)$ of the form $W=R_{a_1}\cdots R_{a_p}$ for some $a_i\in Kan(n)$ such that
\bes
v(I)W&=&v,\\
\overline{v(J)}W&=&0 \hbox{\quad for all } J \subseteq I_n,\\
v(J)W&=&0 \hbox{\quad for all  }  J\subseteq I_n \hbox{ such that } J\neq I \hbox{ as sets}.
\ees
Similarly, there exists $W'=W'(I)$ of the form $W'=R_{b_1}\cdots R_{b_s}$ for some $b_i\in Kan(n)$ such that
\bes
\overline{v(I)}W'&=&v,\\
{v(J)}W'&=&0 \hbox{\quad for all } J \subseteq I_n,\\
\overline{v(J)}W'&=&0 \hbox{\quad for all } J\subseteq I_n \hbox{ such that } J\neq I \hbox{ as sets}.
\ees
\end{lem}
\Proof
Assume, for simplicity, that $I=\{1,\ldots,k\}$. Consider $W_1=R_{e_I}R_{\bar 1}R_{\bar e_{k+1}}\cdots R_{\bar 1}R_{\bar e_n}$; then
we have
\bes
v(I)W_1&=&vR_{\bar 1}R_{\bar e_{k+1}}\cdots R_{\bar 1}R_{\bar e_n}=v(\{ k+1,\ldots,n\}),\\
\overline{v(I)} W_1&=&\pm( v\bar 1)R_{\bar 1}R_{\bar e_{k+1}}\cdots R_{\bar 1}R_{\bar e_n}=\cdots =
\begin{cases}
0  &\hbox{ if } I\neq I_n,\\
\pm v\bar 1&\hbox{ if } I=I_n,
\end{cases}
\\
v(J)W_1&=&
\begin{cases}
0 &\hbox{ if } I\not\subseteq J,\\
{v(J\setminus I)}R_{\bar 1}R_{\bar e_{k+1}}\cdots R_{\bar 1}R_{\bar e_n}& \hbox{ otherwise}
\end{cases}
=\cdots = 0 \hbox{ if } J\neq I,\\
\overline{v(J)}W_1&=&
\begin{cases}
0 &\hbox{ if } I\not\subseteq J,\\
{\pm\a v(J\setminus I)}R_{\bar e_{k+1}}\cdots R_{\bar 1}R_{\bar e_n}& \hbox{ otherwise}
\end{cases}
\\
&=&
\begin{cases}
\pm\a v(\{ k+2,\ldots,n\}) &\hbox{ if } I\subseteq J \hbox{ and } J\setminus I = \{ k+1 \},\\
0 & \hbox{ otherwise}
\end{cases}
\ees
Now, if $I\neq I_n$, we can take $W(I)=W_1R_{e_{\{I_n\setminus I\}}}$, and for $I=I_n$ we can take $W(I)=W_1R_{\bar 1}R_{\bar e_1}R_{e_1}$.
It is easy to check that in both cases $W(I)$ satisfies the conclusions of the first claim of  lemma.

\smallskip

For the second claim, if $\a\neq 0$, it suffices to take $W'(I)=R_{\bar 1}W(I)$. Nevetherless, we will give a general proof.
Assume first that $I\neq I_n$ and let  $i\not\in I$. Consider $\overline {v(I)}\bar e_i=\pm v(I\cup\{i\})$, then 
the element $W'=R_{\bar e_i}W(I\cup\{i\})$, up to sign, satisfies the needed conditions. Finally, for $I=I_n$ one can take 
$W'(I_n)=-R_{e_n}W'(I_n\setminus \{n\})$. 

\ctd

\smallskip
Lemma \ref{lem7} implies several corollaries.

\begin{Cor}
Let $V=V(v,\alpha)$ be the bimodule over $Kan(n)$ with the action defined in  Theorem \ref{thm1}. Then the set of elements (\ref{base}), when $I$ runs all different (non-ordered) subsets of $I_n$, is a base of the  vector space $V$.
\end{Cor}
\Proof
We have already seen  that the set (\ref{base}) generates $V$. Assume that there exists a linear combination
\bes
\sum\limits_{I \subseteq I_n} \beta_I v(I) + \sum\limits_{J \subseteq I_n} \bar\beta_{J} \overline{v(J)} = 0,
\ees
where $\b_I , \bar\beta_J \in F$. Applying the operators $W(I)$ and $W'(J)$, we get that all $\b_I=\bar\beta_J=0$.

\ctd

\begin{Cor}\label{corv}
A special element $v\in V(v,\a)$ is defined uniquely up to a nonzero scalar.
\end{Cor}
\Proof 
Let $v'$ be another special element in $V=V(v,\a)$:
\bes 
v'=\sum\limits_{I \subseteq I_n} \beta_I v(I) + \sum\limits_{J \subseteq I_n} \bar\beta_{J} \overline{v(J)}
\ees
for some $\b_I , \bar\beta_J \in F$. Observe that the operators $W,W'$ in Lemma \ref{lem7} do not depend on choose of the element $v$, and have the same form for elements $v,v'$. Applying the operator $W(\phi)$  to both parts of the above equality, we get  $v'=\beta_{\phi}v$. Clearly, $\beta_{\phi}\neq 0.$

\ctd

\begin{Cor}
The bimodule $V(v,\alpha)$ is irreducible for any $\a\in F$.
\end{Cor}

\Proof
Suppose that  $M$ is a non-zero  sub-bimodule of $V=V(v,\a)$, and choose $0 \ne x \in M$:
\bes
x = \sum\limits_{I \subseteq I_n} \beta_I v(I) + \sum\limits_{J \subseteq I_n} {\bar\beta_J} \overline{v(J)},\ \  \beta_I , \bar\beta_J \in F.
\ees
Since $x \ne 0$, there is some $\beta_I \ne 0$ or $\bar\beta_J \ne 0$. Applying operators $W(I)$ or $W'(J)$ from Lemma \ref{lem7}, we get in both cases that $v\in M$ and hence $M=V$.

\ctd

\begin{rem}
It is easy to check that $V(v,\a)$ for $\a=0$ is isomorphic to the regular bimodule $Reg (Kan(n))$, hence this corollary gives an alternative proof that the superalgebra $Kan(n)$ is simple. 
\end{rem}

\begin{Cor}\label{coraa'}
Bimodules $V(v,\a)$ and $V(v',\a')$ are isomorphic if and only if $|v|=|v'|$ and $\a=\a'$.
\end{Cor}
\Proof
 Denote $V=V(v,\a)$ and $V'=V(v',\a')$. Observe that the operators $W,W'$ in Lemma \ref{lem7} have the same form for elements $v,v'$. Assume  that  $\varphi : V \to V'$ is an isomorphism of bimodules over $Kan(n)$, then we have in $V'$ 
\begin{center}
$\varphi(v) = \sum\limits_{I \subseteq I_n} \beta_I v'(I) + \sum\limits_{I \subseteq I_n} \bar{\b}_I \overline{v'(I)}$, for some $\beta_I , \overline{\beta_I} \in F$.
\end{center}
Applying to both parts of this equality the operator $W=W(\phi)$,  we get $\f(v)=\b_{\phi}v'$, with $0 \ne \beta_{\phi}\in F$. Since $\f$ maintains parity, this is impossible if  $|v|\neq |v'|$. Therefore, $|v|=|v'|$, and we have  
\bes
\a\b_{\phi}v'=\a\f(v)=\f(\a v)=\f(vR_{\bar 1}^2)=\f(v)R_{\bar 1}^2=\a'\f(v)=\a'\b_{\phi}v',
\ees
hence $\a=\a'$.
\ctd

Since, for a given $\a\in F$, the bimodule $V(v,\a)$ is defined, up to isomorphism, by the parity of $v$, we will denote by $V(\a)$ the bimodule $V(v,\a)$ with $|v|=|n|$. 

\smallskip
Recall that, for a superalgebra $A=A_0\oplus A_1$, an $A$-superbimodule
$V^{\mathrm{op}} = V_{ 0}^{\mathrm{op}} + V_{ 1}^{\mathrm{op}}$ is called   {\it opposite} to an $A$-superbimodule $V=V_0\oplus V_1$, if
$V_{0}^{\mathrm{op}} = V_{ 1}$,
$V_{ 1}^{\mathrm{op}} = V_{ 0}$
and $A$ acts on it by the following rule:
$a\cdot v = (-1)^{|a|} av$, $v\cdot a = va$,
where $v \in V^{\mathrm{op}},\ a \in A_{0}\cup A_{1}$.

It is easy to check that, for any  superbimodule~$V$, the identical application $V\rightarrow V^{op},\ v\mapsto v$, is an odd isomorphism between $V$ and $V^{op}$. In particular, if $V$ is Jordan, the opposite superbimodule~$V^{\mathrm{op}}$  is Jordan as well. We sometimes will say that the bimodule $V^{op}$ is obtained from $V$ by {\it changing of parity}.

Corollaries \ref{corv} and \ref{coraa'} imply
\begin{prop}\label{prop1}
Every unital irreducible bimodule over $Kan(n)$ is isomorphic to a bimodule $V(\a)$ or to its opposite $V(\a)^{\mathrm{op}} $. Moreover, the bimodules $V(\a)$ and $V(\a)^{\mathrm{op}}$ are not isomorphic. 
\end{prop}
\Proof It suffices to note that, for $|v'|=|v|+1$, the mapping 
\bes
\sum\limits_{I \subseteq I_n} \beta_I v(I) + \sum\limits_{J \subseteq I_n} {\bar\beta_J} \overline{v(J)}\mapsto \sum\limits_{I \subseteq I_n} \beta_I v'(I) + \sum\limits_{J \subseteq I_n} {\bar\beta_J} \overline{v'(J)}
\ees
defines an isomorphism of the bimodules $V(v,\a)^{op}$ and $V(v',\a)$.

\ctd

\section{$V(\a)$ is Jordan}

\hspace{\parindent}
Finally, we will show that $V(\a)$ is a Jordan bimodule over $Kan(n)$ for all $\alpha$. For this, we will  embedd it into a Jordan superalgebra.

\smallskip

Recall than a linear operator $E$ on a unital algebra $A$ is called  {\em a generalized derivation} of $A$ if it satisfies the relation
\bes
E(ab)=E(a)b+aE(b)-abE(1).
\ees
Let $P=\la P_0\oplus P_1,\cdot,\{,\}\ra$ be a  unital Poisson superalgebra, $E:P\rightarrow P$  be a generalized derivation of $P$ which satisfies also the condition
\bee\label{Pder}
E(\{p,q\})=\{E(p),q\}+\{p,E(q)\}+\{p,q\}E(1).
\eee
Furthermore, let $(A,D)$ be a commutative associative  algebra with a derivation $D$. Define the following bracket on the tensor product $P\otimes A$:
\bee\label{br}
\la p\otimes a,q\otimes b\ra&=& \{p,q\}\otimes ab+E(p)q\otimes aD(b)-(-1)^{|p||q|}E(q)p\otimes D(a)b,
\eee
where $p,q\in P,\ a,b\in A$. 
\begin{thm}\label{bracket}
The bracket (\ref{br}) is a Jordan bracket on the commutative and associative superalgebra $P\otimes A=(P_0\otimes A)\oplus (P_1\otimes A)$.
\end{thm}
\Proof 
Observe first that a commutative associative superalgebra $P\otimes A$ with a
superanticommutative bracket $<,>$ satisfies graded identities (2)-(4) if and only if the Grassmann envelope $G(P\otimes A)=G_0\otimes (P_0\otimes A)+G_1\otimes (P_1\otimes A)$, with the
bracket $\la a\otimes g, b\otimes h\ra=\la a,b\ra\otimes gh$, satisfies the
nongraded versions of these identities. It is easy to check the isomorphism
\bes
G(P\otimes A)\cong G(P)\otimes A,
\ees
where $G(P)=G_0\otimes P_0+G_1\otimes P_1$ is the Grassmann envelope of the superalgebra $P$. So, passing to the Grassmann
envelope, we see that it sufficient to prove nongraded identities (2)--(4) for the case  when $P$ is a Poisson algebra (not a superalgebra).

Let us first check  identity (\ref{gleibniz}):
\bes
&\la (p\otimes a)(q\otimes b),r\otimes c\ra=
\la pq\otimes ab,r\otimes c\ra&\\
&= \{pq,r\}\otimes abc+E(pq)r\otimes abD(c)-E(r)pq\otimes cD(ab)&\\
&=(p\{q,r\}+q\{p,r\})\otimes abc+(E(p)q+pE(q)-pqE(1))r\otimes abD(c)-E(r)pq\otimes cD(ab)&\\
&=(p\otimes a)(\{q,r\}\otimes  bc +E(q)r\otimes bD(c)-qE(r)\otimes cD(b))&\\
&+(q\otimes b)(  \{p,r\}\otimes  ac +E(p)r\otimes aD(c)-pE(r)\otimes cD(a))-pqE(1)r\otimes abD(c)&\\
&=(p\otimes a)\la q\otimes b,r\otimes c\ra+ (q\otimes b)\la p\otimes a,r\otimes c\ra
-(p\otimes a) (q\otimes b)\la 1,r\otimes c\ra,&
\ees
so (\ref{gleibniz}) is satisfied.

 Furthermore, the nongraded version of identity (\ref{gjacob}) has form
\bes
J(a,b,c)=S(a,b,c),
\ees
where $J(a,b,c)=\la\la a,b\ra,c\ra + \la\la b,c\ra,a\ra+\la\la c,a\ra,b\ra$ and $S(a,b,c)=\la a,b\ra \la 1,c\ra +\la b,c\ra \la 1,a\ra+\la c,a\ra \la 1,b\ra$. 

Consider  
\bes
&\la \la p\otimes a,q\otimes b\ra,r\otimes c\ra=\la \{p,q\}\otimes ab+E(p)q\otimes aD(b)-E(q)p\otimes bD(a),r\otimes c\ra&\\
&=\{\{p,q\},r\}\otimes abc+E(\{p,q\})r\otimes abD(c)-E(r)\{p,q\}\otimes cD(ab)\\
&+\{E(p)q,r\}\otimes aD(b)c+E(E(p)q)r\otimes aD(b)D(c)-E(r)E(p)q\otimes D(aD(b))c&\\
&-\{E(q)p,r\}\otimes bD(a)c-E(E(q)p)r\otimes bD(a)D(c)+E(r)E(q)p\otimes D(bD(a))c.&
\ees
By  the properties of  the bracket $\{,\}$ and the generalized derivation $E$, we have further
\bes
&\la \la p\otimes a,q\otimes b\ra,r\otimes c\ra=\{\{p,q\},r\}\otimes abc+(\{E(p),q\}+\{p,E(q)\}+\{p,q\}E(1))r\otimes abD(c)&\\
&-E(r)\{p,q\}\otimes c(D(a)b+aD(b))+(E(p)\{q,r\}+q\{E(p),r\})\otimes aD(b)c&\\
&+(E^2(p)q+E(p)E(q)-E(p)qE(1))r\otimes aD(b)D(c)&\\
&-E(r)E(p)q\otimes (D(a)D(b)c+aD^2(b)c)-(E(q)\{p,r\}-p\{E(q),r\})\otimes bD(a)c&\\
&-(E^2(q)p+E(q)E(p)-E(q)pE(1))r\otimes bD(a)D(c)&\\
&+E(r)E(q)p\otimes (D(b)D(a)+bD^2(a))c&.
\ees
Calculating the cyclic sum, we get
\bes
J(p\otimes a,q\otimes b,r\otimes c)=(E(1)\otimes 1)(\{p,q\}r\otimes abD(c)+
\{q,r\}p\otimes bcD(a)+\{r,p\}q\otimes caD(b)).
\ees
Consider now
\bes 
&\la p\otimes a,q\otimes b\ra \la 1\otimes 1,r\otimes c\ra=(\{p,q\}\otimes ab+E(p)q\otimes aD(b)-E(q)p\otimes D(a)b)(E(1)r\otimes D(c))&\\
&=(E(1)\otimes 1)(\{p,q\}r\otimes abD(c)+E(p)qr\otimes aD(b)D(c)-E(q)pr\otimes D(a)bD(c)).&
\ees
Therefore, the cyclic sum
\bes
S(p\otimes a,q\otimes b,r\otimes c)&=&(E(1)\otimes 1)(\{p,q\}r\otimes abD(c)+
\{q,r\}p\otimes bcD(a)+\{r,p\}q\otimes caD(b))\\
&=&J(p\otimes a,q\otimes b,r\otimes c),
\ees
proving (\ref{gjacob}). 

Finally, observe that all partial linearizations of (4) follow from 
 identity (\ref{gjacob}), hence to prove (4) it suffices for us to prove that
for any $p\in P_1$ and $a\in A$ holds
\bes
\la\la p\otimes a,p\otimes a\ra,p\otimes a\ra=\la p\otimes a,p\otimes a\ra\la 1\otimes 1,p\otimes a\ra.
\ees
We have
\bes
\la p\otimes a,p\otimes a\ra&=& \{p,p\}\otimes aa+E(p)p\otimes aD(a)+E(p)p\otimes D(a)a\\
&=&\{p,p\}\otimes a^2+E(p)p\otimes D(a^2),
\ees
and, furthermore,
\bes
&\la\la p\otimes a,p\otimes a\ra,p\otimes a\ra=\la \{p,p\}\otimes a^2+E(p)p\otimes D(a^2),p\otimes a\ra&\\
&=\{\{p,p\},p\}\otimes a^3+E(\{p,p\})p\otimes a^2D(a)-E(p)\{p,p\}\otimes D(a^2)a&\\
&+\{E(p)p,p\}\otimes D(a^2)a+E(E(p)p)p\otimes D(a^2)D(a)-E(p)E(p)p\otimes D^2(a^2)a.&
\ees
We have $\{\{p,p\},p\}=p^2=E(p)^2=0$, therefore by (\ref{Pder})
\bes
&\la\la p\otimes a,p\otimes a\ra,p\otimes a\ra=(2\{E(p),p\}+E(1)\{p,p\})p\otimes a^2D(a)-E(p)\{p,p\}\otimes D(a^2)a&\\
&+(E(p)\{p,p\}-p\{E(p),p\})\otimes D(a^2)a=E(1)\{p,p\}p\otimes a^2D(a).&
\ees
On the other hand,
\bes 
\la p\otimes a,p\otimes a\ra \la 1\otimes 1,p\otimes a\ra
&=&
(\{p,p\}\otimes a^2+E(p)p\otimes D(a^2))(E(1)p\otimes D(a))\\
&=&E(1)\{p,p\}p\otimes a^2D(a).
\ees
Hence (4) is true, and the theorem is proved.

\ctd

\begin{Cor}\label{gr_Poisson}
Let $P=\oplus_{i=0}^{\infty}P_i$ be a $\Z$-graded Poisson superalgebra such that $\{P_i,P_j\}\subseteq P_{i+j-2}$.
Then the application $E:P\rightarrow P,\ E:p_i\mapsto (i-1)p_i, \ p_i\in P_i,$ is a generalized derivation of $P$ which satisfies relation (\ref{Pder}).
In particular, for any associative and commutative algebra $(A,D)$ with a derivation $D$, the tensor product superalgebra $P\otimes A$ has a Jordan bracket given by (\ref{br}).
\end{Cor}

The Grassmann superalgebra $G_n$ has a natural $\Z$-grading given by degrees of monomials:  $e_I\in(G_n)_i$ if and only if $|I|=i$. Clearly, $\{(G_n)_i,(G_n)_j\}\subseteq (G_n)_{i+j-2}$, hence $G_n$ satisfies the previous corollary.
Consider the polynomial algebra $A=F[t]$ with a derivation $D_{\a}$ defined by the condition $D_{\a}(t)=-\a t$, then  the superalgebra $G_n[t]\cong G_n\otimes F[t]$ has a Jordan bracket defined by  (\ref{br}) with $D=D_{\a}$. Therefore, we have

\begin{Cor}\label{cor_grassmann}
The Kantor double $J(G_n[t])_{\a}$ with respect to the bracket defined on $G_n[t]=G_n\otimes F[t]$ according to (\ref{br}) with the derivation $D_{\a}$, is a Jordan superalgebra.  
\end{Cor}

\smallskip

Now, we will find in the superalgebra $J(G_n[t])_{\a}$ a $Kan(n)$-subbimodule isomorphic to the bimodule $V(\a)$. Since $\la a,b\ra=\{a,b\}$ for $a,b\in G_n$,
the superalgebra $Kan(n)=J(G_n)$ is a subsuperalgebra of $J(G_n[t])_{\a}$.
Consider in $J(G_n[t])_{\a}$ the subspace
$W=G_n\otimes t+\overline{G_n\otimes t}$. Clearly, $W\bullet G_n\subseteq G_n$, and for $a,b\in G_n$ we have
\bes 
\overline{ a\otimes 1}\bullet  \overline{b\otimes t}&=&(-1)^{|b|}\la a\otimes 1,b\otimes t\ra=(-1)^{|b|}(\{a,b\}\otimes t+E(a)b\otimes \a t)\\
&=&(-1)^{|b|}\a\{a,b\}E(a)b\otimes t\in W.
\ees
Therefore, $W$ is a unital Jordan bimodule over $Kan(n)=J(G_n)$. Let $w=e_{I_n}\otimes t$, then it is clear that $w$ is a specal element in $W$ that satisfies the properties of Lemma \ref{lem3}.
Moreover, one can easily check that $W=V(w,\a)$, in notation of Section 5; the exact isomorphism is given by the mapping
\bes 
v(I)\mapsto (-1)^{sgn\,\s}e_{\{I_n\setminus I\}}\otimes t,\ \overline{v(I)}\mapsto (-1)^{sgn\,\s}\overline{e_{\{I_n\setminus I\}}\otimes t},
\ees
where, for $I=\{i_1,\ldots,i_k\}$,  $\s$ is a permutation $\s:I_n\mapsto (I_n\setminus I)\cup\{i_k,\ldots,i_1\}$.

\smallskip

Resuming,  we can formulate our main theorem:

\begin{thm}
The bimodule $V_{\a}$ is a unital Jordan irreducible bimodule over the superalgebra $Kan(n)$, and every such a bimodule over $Kan(n)$ for $n\geq 2$ over an algebraically closed field of characteristic not 2  is isomorphic to $V_{\a}$ or to its opposite bimodule.
\end{thm}

\section{Acknowledgements}
The paper is a part of the PhD thesis of the first author done at the University of S\~ao Paulo. He acknowledges the support by the  CAPES and CNPq grants.
The second author acknowledges the support by the FAPESP and CNPq grants.

\end{document}